\documentclass[11pt,thmsa,arabtex]{article}
\usepackage{amssymb}
\usepackage{amsfonts}
\usepackage{amsmath}
\usepackage{graphicx}
\usepackage{hyperref}
\usepackage{cite}
\usepackage{ntheorem}
\usepackage{epstopdf}

\newtheorem*{theorem-non}{Theorem}

\newtheorem{definition}{Definition}[section]

\newtheorem{example}{Example}[section]

\numberwithin{theorem}{section}
\numberwithin{equation}{section}

\begin{document}

\title{{\Large \textbf{Smarandache curves of some special curves in the Galilean 3-space}}}
\author{\emph{\textbf{H. S. Abdel-Aziz}} and \emph{\textbf{M. Khalifa Saad}}\thanks{
~E-mail address:~mohamed\_khalifa77@science.sohag.edu.eg} \\
{\small \emph{Dept. of Math., Faculty of Science, Sohag Univ., 82524 Sohag,
Egypt}}}
\date{}
\maketitle

\textbf{Abstract.} In the present paper, we consider a position vector of an arbitrary curve in the three-dimensional Galilean space $G_{3}$. Furthermore, we give some conditions on the curvatures of this arbitrary curve to study special curves and their Smarandache curves. Finally, in the light of this study, some related examples of these curves
are provided and plotted.

\section{Introduction}

Discovering Galilean space-time is probably one of the major achievements of non relativistic physics. Nowadays Galilean space is becoming increasingly popular as evidenced from the connection of the fundamental concepts such as velocity, momentum, kinetic energy, etc. and principles as indicated in \cite{IY}. In recent years, researchers have begun to investigate curves and surfaces in the Galilean space and thereafter pseudo-Galilean space.\\
 In classical curve theory, the geometry of a curve in three-dimensions is essentially characterized by two scalar functions, curvature $\kappa$ and torsion $\tau$ as well as its Frenet vectors.
 A regular curve in Euclidean space whose position vector is composed by Frenet frame vectors on another regular curve is called a Smarandache curve. Smarandache curves have been investigated by some differential geometers
\cite{AT,TY}. M. Turgut and S. Yilmaz  defined a special case of such curves and call it Smarandache $TB_{2}$ curves in the space $E_{1}^{4}$ \cite{TY}. They studied special Smarandache curves which are defined by the tangent and second binormal vector fields. Additionally, they computed formulas of this kind curves. In \cite{AT}, the author introduced some special Smarandache curves in the Euclidean space. He studied Frenet-Serret invariants of a special case.\\
  In the field of computer aided design and computer graphics, helices can be used for the tool path description, the simulation of kinematic motion or the design of highways, etc. \cite{XY}. The main feature of general helix or slope line is that the tangent makes a constant angle with a fixed direction in every point which is called the axis of the general helix. A classical result stated by Lancret in $1802$ and first proved by de Saint Venant in 1845 says that: A necessary and sufficient condition that a curve be a general helix is that the ratio ($\kappa / \tau$) is constant along the curve, where $\kappa$ and $\tau$ denote the curvature and the torsion, respectively. Also, the helix is also known as circular helix or W-curve which is a special case of the general helix \cite{DS}.\\
  Salkowski (resp. Anti-Salkowski) curves in Euclidean space are generally known as family of curves with constant curvature (resp. torsion) but nonconstant torsion (resp. curvature) with an explicit parametrization.They were defined in an earlier paper \cite{ES}.\\
  In this paper, we compute Smarandache curves for a position vector of an arbitrary curve and some of its special curves. Besides, according to Frenet frame {$\textbf{T}$, $\textbf{N}$, $\textbf{B}$} of the considered curves in the Galilean space $G_{3}$, the meant Smarandache curves $\textbf{TN}$, $\textbf{TB}$ and $\textbf{TNB}$ are obtained.
   We hope these results will be helpful to mathematicians who are specialized on mathematical modeling.
\section{Preliminaries}
Let us recall the basic facts about the three-dimensional Galilean geometry $%
G_{3}$. The geometry of the Galilean space has been firstly explained in
\cite{OR}. The curves and some special surfaces in $G_{3}$ are considered in
\cite{DM}. The Galilean geometry is a real Cayley-Klein geometry with
projective signature $(0,0,+,+)$ according to \cite{EM}. The absolute of the
Galilean geometry is an ordered triple $({w,f,I}$ $)$ where $w$ is the ideal
(absolute) plane $(x_{0}=0)$, $f$ is a line in $w$ $(x_{0}=x_{1}=0)$ and $I$
is elliptic $((0:0:x_{2}:x_{3})\longrightarrow (0:0:x_{3}:-x_{2}))$
involution of the points of $f$ . In the Galilean space there are just two
types of vectors, non-isotropic $\mathbf{x}(x,y,z)$ (for which holds $x\neq
0 $). Otherwise, it is called isotropic. We do not distinguish classes of
vectors among isotropic vectors in $G_{3}$. A plane of the form $x=const$.
in the Galilean space is called Euclidean, since its induced geometry is
Euclidean. Otherwise it is called isotropic plane. In affine coordinates,
the Galilean inner product between two vectors $P=(p_{1},p_{2},p_{3})$ and $%
Q=(q_{1},q_{2},q_{3})$ is defined by \cite{PK}:
\begin{equation}
\langle P,Q\rangle _{G_{3}}=\left\{
\begin{array}{c}
p_{1}q_{1}\ \ \ \ \ \ \ \ \ \ \ \ \text{if}\ p_{1}\neq 0\vee q_{1}\neq 0, \\
p_{2}q_{2}+p_{3}q_{3}\ \ \text{if\ }p_{1}=0\wedge q_{1}=0.%
\end{array}%
\right.
\end{equation}%
And the cross product in the sense of Galilean space is given by:
\begin{equation}
\left( P\times Q\right) _{G_{3}}=\left\{
\begin{array}{c}
\left\vert
\begin{array}{ccc}
0 & e_{2} & e_{3} \\
p_{1} & p_{2} & p_{3} \\
q_{1} & q_{2} & q_{3}%
\end{array}%
\right\vert \ ;\ \ \ \ \ \text{if}\ p_{1}\neq 0\vee q_{1}\neq 0, \\
\\
\left\vert
\begin{array}{ccc}
e_{1} & e_{2} & e_{3} \\
p_{1} & p_{2} & p_{3} \\
q_{1} & q_{2} & q_{3}%
\end{array}%
\right\vert \ ;\ \ \ \ \ \text{if\ }p_{1}=0\wedge q_{1}=0.
\end{array}
\right.
\end{equation}
A curve $\eta (t)=(x(t),y(t),z(t))$ is admissible in $G_{3}$ if it has no
inflection points $(\dot{\eta}(t)\times \ddot{\eta}(t)\neq 0)$ and no isotropic tangents $(\dot{x}(t)\neq 0)$. An
admissible curve in $G_{3}$ is an analogue of a regular curve in Euclidean
space.
For an admissible curve $\eta $ $:I$ \ $\rightarrow $\ $G_{3},$ $I\subset
R $ parameterized by the arc length $s$ with differential form $dt=ds$,
given by
\begin{equation}
\eta (s)=(s,y(s),z(s)).
\end{equation}
The curvature $\kappa (s)$ and torsion $\tau (s)$ of $\eta $ are defined
by
\begin{eqnarray}
\kappa (s) &=&\left\Vert \eta ^{^{\prime \prime }}(s)\right\Vert =\sqrt{%
y^{^{\prime \prime }}(s)^{2}+z^{^{\prime \prime }}(s)^{2}},  \notag \\
\tau (s) &=&\frac{det(\eta^{\prime}(s),\eta^{\prime \prime}(s),\eta^{\prime \prime \prime}(s))}{\kappa ^{2}(s)
}.
\end{eqnarray}
Note that an admissible curve has non-zero curvature.
The associated trihedron is given by
\begin{eqnarray}
\mathbf{T}(s) &=&\eta ^{^{\prime }}(s)=(1,y^{^{\prime }}(s),z^{^{\prime
}}(s)),  \notag \\
\mathbf{N}(s) &=&\frac{\eta ^{^{\prime \prime }}(s)}{\kappa (s)}=\frac{%
(0,y^{^{\prime \prime }}(s),z^{^{\prime \prime }}(s))}{\kappa (s)},  \notag
\\
\mathbf{B}(s) &=&\frac{(0,-z^{^{\prime \prime }}(s),y^{^{\prime \prime }}(s))%
}{\kappa (s)}.
\end{eqnarray}
For derivatives of the tangent $\mathbf{T}$, normal $\mathbf{N}$ and
binormal $\mathbf{B}$ vector field, the following Frenet formulas in the
Galilean space hold \cite{OR}
\begin{equation}
\left[
\begin{array}{c}
\mathbf{T} \\
\mathbf{N} \\
\mathbf{B}%
\end{array}%
\right] ^{{\large {\prime }}}=\left[
\begin{array}{ccc}
0 & \kappa & 0 \\
0 & 0 & \tau \\
0 & -\tau & 0%
\end{array}%
\right] \left[
\begin{array}{c}
\mathbf{T} \\
\mathbf{N} \\
\mathbf{B}%
\end{array}%
\right] .
\end{equation}%
From $(2.5)$ and $(2.6)$, we derive an important relation
\begin{equation*}
\eta ^{\prime \prime \prime }(s)=\kappa ^{\prime }(s)\mathbf{N}(s)+\kappa
(s)\tau (s)\mathbf{B}(s).
\end{equation*}
In \cite{TY} authors introduced:
\begin{definition}
A regular curve in Minkowski space-time, whose position vector is composed by Frenet frame vectors on another regular curve, is called a Smarandache curve.
\end{definition}
In the light of the above definition, we adapt it to admissible curves in the Galilean space as
follows:
\begin{definition}
let $\eta =\eta (s)$ be an admissible curve in $G_{3}$ and $\{\mathbf{T},
\mathbf{N},\mathbf{B}\}$ be its moving Frenet frame. Smarandache $\mathbf{TN}
,\mathbf{TB}$ and $\mathbf{TNB}$ curves are respectively, defined by\\
\begin{eqnarray}
\eta _{\mathbf{TN}} &=&\frac{\mathbf{T}+\mathbf{N}}{\left\Vert \mathbf{T}+%
\mathbf{N}\right\Vert },\ \ \   \notag \\
\eta _{\mathbf{TB}} &=&\frac{\mathbf{T}+\mathbf{B}}{\left\Vert \mathbf{T}+%
\mathbf{B}\right\Vert },  \notag \\
\eta _{\mathbf{TNB}} &=&\frac{\mathbf{T}+\mathbf{N+B}}{\left\Vert \mathbf{T%
}+\mathbf{N+B}\right\Vert }.
\end{eqnarray}
\end{definition}

\section{Smarandache curves of an arbitrary curve in $G_{3}$}

In this section, we consider the position vector of an arbitrary curve with
curvature $\kappa (s)$ and torsion $\tau (s)$ in the Galilean space $G_{3}$
which introduced by \cite{AT} as follows
\begin{equation}
\textbf{r}(s)=\left(
s,\int \left( \int \kappa (s)\cos \left( \int \tau (s)ds\right) ds\right) ds,
\int \left( \int \kappa (s)\sin \left( \int \tau (s)ds\right) ds\right) ds
\right).
\end{equation}
The derivatives of this curve are respectively, given by
\begin{equation*}
\textbf{r}^{\prime }(s)=\left(
1,\int \kappa (s)\,\cos \left( \int \tau (s)\,ds\right) ds,
\int \kappa (s)\,\sin \left( \int \tau (s)\,ds\right) \,ds
\right),
\end{equation*}
\begin{equation*}
\textbf{r}^{\prime \prime }(s)=\left( 0,\kappa (s)\cos \left( \int \tau
(s)\,ds\right) ,\kappa (s)\,\sin \left( \int \tau (s)\,ds\right) \right),
\end{equation*}
\begin{equation}
\textbf{r}^{\prime \prime \prime }(s)=\left(
\begin{array}{c}
0,\kappa ^{\prime }(s)\cos \left( \int \tau (s)\,ds\right) -\kappa (s)\tau
(s)\sin \left( \int \tau (s)\,ds\right) , \\
\kappa ^{\prime }(s)\sin \left( \int \tau (s)\,ds\right) +\kappa (s)\tau
(s)\cos \left( \int \tau (s)\,ds\right)%
\end{array}
\right).
\end{equation}
The frame vector fields of $\textbf{r}$ are as follows\\

\begin{equation*}
\textbf{T}_{\textbf{r}}=\left(
1,\int \kappa (s)\cos \left( \int \tau (s)\,ds\right) ds,
\int \kappa (s)\sin \left( \int \tau (s)\,ds\right) \,ds
\right),
\end{equation*}
\begin{equation*}
\textbf{N}_{\textbf{r}}=\left( 0,\cos \left( \int \tau (s)\,ds\right) ,\sin \left( \int \tau
(s)\,ds\right) \right),
\end{equation*}
\begin{equation}
\textbf{B}_{\textbf{r}}=\left( 0,-\sin \left( \int \tau (s)\,ds\right) ,\cos \left( \int \tau
(s)\,ds\right) \right).
\end{equation}
By Definition (2.2), the $\textbf{TN}$, $\textbf{TB}$ and $\textbf{TNB}$ Smarandache curves of $r$ are respectively, written as
\begin{equation*}
\textbf{r}_{\textbf{TN}}=\left(
\begin{array}{c}
1,\cos \left( \int \tau (s)\,ds\right) +\int \kappa (s)\cos \left( \int \tau
(s)\,ds\right) \,ds, \\
\int \kappa (s)\sin \left( \int \tau (s)\,ds\right) \,ds+\sin \left( \int
\tau (s)\,ds\right)%
\end{array}%
\right),
\end{equation*}
\begin{equation*}
\textbf{r}_{\textbf{TB}}=\left(
\begin{array}{c}
1,\int \kappa (s)\cos \left( \int \tau (s)\,ds\right) \,ds-\sin \left( \int
\tau (s)\,ds\right) , \\
\cos \left( \int \tau (s)\,ds\right) +\int \kappa (s)\sin \left( \int \tau
(s)\,ds\right) ds%
\end{array}%
\right),
\end{equation*}
\begin{equation}
\textbf{r}_{\textbf{TNB}}=\left(
\begin{array}{c}
1,\cos \left( \int \tau (s)\,ds\right) +\int \kappa (s)\cos \left( \int \tau
(s)\,ds\right) \,ds \\
-\sin \left( \int \tau (s)\,ds\right) ,\cos \left( \int \tau (s)\,ds\right) +
\\
\int \kappa (s)\sin \left( \int \tau (s)\,ds\right) ds+\sin \left( \int \tau
(s)\,ds\right)
\end{array}%
\right).
\end{equation}

\section{Smarandache curves of some special curves in $G_{3}$}
\subsection{Smarandache curves of a general helix}

Let $\alpha (s)$ be a general helix in $G_{3}$ with ($\tau /\kappa
=m=const.$)which can be written as
\begin{equation}
\alpha (s)=\left(
\begin{array}{c}
s,\frac{1}{m}\int \sin \left( m\int \kappa (s)\,ds\right) \,ds, \\
\frac{-1}{m}\int \cos \left( m\int \kappa (s)\,ds\right) \,ds
\end{array}
\right).
\end{equation}
Then $\alpha ^{\prime }, \alpha ^{\prime \prime }, \alpha ^{\prime \prime \prime }$ for this curve are respectively, expressed as
\begin{equation*}
\alpha ^{\prime }(s)=\left( 1,\frac{1}{m}\sin \left( m\int \kappa
(s)\,ds\right) ,\frac{-1}{m}\cos \left( m\int \kappa (s)\,ds\right) \right),
\end{equation*}
\begin{equation*}
\alpha ^{\prime \prime }(s)=\left( 0,\kappa (s)\cos \left( m\int \kappa
(s)\,ds\right) ,\kappa (s)\sin \left( m\int \kappa (s)\,ds\right) \right),
\end{equation*}
\begin{equation}
\alpha ^{\prime \prime \prime }(s)=\left(
\begin{array}{c}
0,\kappa ^{\prime }(s)\cos \left( m\int \kappa (s)\,ds\right) - \\
m~\kappa ^{2}(s)\sin \left( m\int \kappa (s)\,ds\right) , \\
\kappa ^{\prime }(s)\sin \left( m\int \kappa (s)\,ds\right) + \\
m~\kappa ^{2}(s)\cos \left( m\int \kappa (s)\,ds\right)
\end{array}
\right).
\end{equation}
The moving Frenet vectors of $\alpha (s)$ are given by
\begin{equation*}
\textbf{T}_{\alpha}=\left( 1,\frac{1}{m}\sin \left( m\int \kappa (s)\,ds\right) ,\frac{-1}{m}%
\cos \left( m\int \kappa (s)\,ds\right) \right),
\end{equation*}
\begin{equation*}
\textbf{N}_{\alpha}=\left( 0,\cos \left( m\int \kappa (s)\,ds\right) ,\sin \left( m\int \kappa
(s)\,ds\right) \right),
\end{equation*}
\begin{equation}
\textbf{B}_{\alpha}=\left( 0,-\sin \left( m\int \kappa (s)\,ds\right) ,\cos \left( m\int
\kappa (s)\,ds\right) \right).
\end{equation}
From which, Smarandache curves are obtained
\begin{equation*}
\alpha_{\textbf{TN}}=\left(
\begin{array}{c}
1,\cos \left( m\int \kappa (s)\,ds\right) +\frac{1}{m}\sin \left( m\int
\kappa (s)\,ds\right), \\
\frac{-1}{m}\cos \left( m\int \kappa (s)\,ds\right) +\sin \left( m\int
\kappa (s)\,ds\right)
\end{array}
\right),
\end{equation*}
\begin{equation*}
\alpha_{\textbf{TB}}=\left(
1,-\left( \frac{m-1}{m}\right) \sin \left( m\int \kappa (s)\,ds\right) ,
\left( \frac{m-1}{m}\right) \cos \left( m\int \kappa (s)\,ds\right)
\right),
\end{equation*}
\begin{equation}
\alpha_{\textbf{TNB}}=\left(
\begin{array}{c}
1,\cos \left( m\int \kappa (s)\,ds\right) -\left( \frac{m-1}{m}\right) \sin
\left( m\int \kappa (s)\,ds\right), \\
\left( \frac{m-1}{m}\right) \cos \left( m\int \kappa (s)\,ds\right) +\sin
\left( m\int \kappa (s)\,ds\right)
\end{array}
\right).
\end{equation}

\subsection{Smarandache curves of a circular helix}

Let $\beta (s)$ be a circular helix in $G_{3}$ with ($\tau=a=const., \kappa
=b=const.$) which can be written as
\begin{equation}
\beta (s)=\left( s,a\int \left( \int \cos (bs)\,ds\right) \,ds,a\int \left(
\int \sin (bs)\,ds\right) \,ds\right).
\end{equation}
For this curve, we have
\begin{equation*}
\beta ^{\prime }(s)=\left( 1,\frac{a}{b}\sin (bs),-\frac{a}{b}\cos
(bs)\right),
\end{equation*}
\begin{equation*}
\beta ^{\prime \prime }(s)=\left( 0,a\cos (bs),a\sin (bs)\right),
\end{equation*}
\begin{equation}
\beta ^{\prime \prime \prime }(s)=\left( 0,-ab\sin (bs),ab\cos (bs)\right).
\end{equation}
Making necessary calculations from above, we have
\begin{equation*}
\textbf{T}_{\beta}=\left( 1,\frac{a}{b}\sin (bs),-\frac{a}{b}\cos (bs)\right),
\end{equation*}
\begin{equation*}
\textbf{N}_{\beta}=\left( 0,\cos (bs),\sin (bs)\right),
\end{equation*}
\begin{equation}
\textbf{B}_{\beta}=\left( 0,-\sin (bs),\cos (bs)\right).
\end{equation}
Considering the last Frenet vectors, the $\textbf{TN}$, $\textbf{TB}$ and $\textbf{TNB}$ Smarandache curves of $\beta$ are respectively, as follows
\begin{equation*}
\beta_{\textbf{TN}}=\left( 1,\cos (bs)+\frac{a}{b}\sin (bs),-\frac{a}{b}\cos (bs)+\sin
(bs)\right),
\end{equation*}
\begin{equation*}
\beta_{\textbf{TB}}=\left( 1,\left( \frac{a-b}{b}\right) \sin (bs),\left( \frac{b-a}{b}
\right) \cos (bs)\right),
\end{equation*}
\begin{equation}
\beta_{\textbf{TNB}}=\left(
\begin{array}{c}
1,\cos (bs)+\left( \frac{a-b}{b}\right) \sin (bs), \\
\left( \frac{b-a}{b}\right) \cos (bs)+\sin (bs)
\end{array}
\right).
\end{equation}

\subsection{Smarandache curves of a Salkowski curve}

Let $\gamma (s)$ be a Salkowski curve in $G_{3}$ with ($\tau=\tau (s), \kappa
=a=const.$) which can be written as
\begin{equation}
\gamma (s)=\left(
\begin{array}{c}
s,a\int \left( \int \cos \left( \int \tau (s)\,ds\right) ds\,\right) ds, \\
a\int \left( \int \sin \left( \int \tau (s)\,ds\right) ds\,\right) ds%
\end{array}
\right).
\end{equation}
If we differentiate this equation three times, one can obtain
\begin{equation*}
\gamma ^{\prime }(s)=\left( 1,a\int \cos \left( \int \tau (s)\,ds\right)
ds\,,a\int \sin \left( \int \tau (s)\,ds\right) ds\right),
\end{equation*}
\begin{equation*}
\gamma ^{\prime \prime }(s)=\left( 0,a\cos \left( \int \tau (s)\,ds\right)
,a\sin \left( \int \tau (s)\,ds\right) \right),
\end{equation*}
\begin{equation}
\gamma ^{\prime \prime \prime }(s)=\left(
\begin{array}{c}
0,-a~\tau (s)\sin \left( \int \tau (s)\,ds\right) , \\
a~\tau (s)\cos \left( \int \tau (s)\,ds\right)
\end{array}
\right).
\end{equation}
In addition to that, the tangent, principal normal and binormal vectors of $\gamma$ are in the following forms
\begin{equation*}
\textbf{T}_{\gamma}=\left( 1,a\int \cos \left( \int \tau (s)\,ds\right) ds,a\int \sin \left(
\int \tau (s)\,ds\right) ds\right),
\end{equation*}
\begin{equation*}
\textbf{N}_{\gamma}=\left( 0,\cos \left( \int \tau (s)\,ds\right) ,\sin \left( \int \tau
(s)\,ds\right) \right),
\end{equation*}
\begin{equation}
\textbf{B}_{\gamma}=\left( 0,-\sin \left( \int \tau (s)\,ds\right) ,\cos \left( \int \tau
(s)\,ds\right) \right).
\end{equation}
Furthermore, Smarandache curves for $\gamma$ are
\begin{equation*}
\gamma_{\textbf{TN}}=\left(
\begin{array}{c}
1,\cos \left( \int \tau (s)\,ds\right) +a\int \cos \left( \int \tau
(s)\,ds\right) ds, \\
a\int \sin \left( \int \tau (s)\,ds\right) ds+\sin \left( \int \tau
(s)\,ds\right)
\end{array}
\right),
\end{equation*}
\begin{equation*}
\gamma_{\textbf{TB}}=\left(
\begin{array}{c}
1,a\int \cos \left( \int \tau (s)\,ds\right) ds-\sin \left( \int \tau
(s)\,ds\right) , \\
\cos \left( \int \tau (s)\,ds\right) +a\int \sin \left( \int \tau
(s)\,ds\right) ds
\end{array}
\right),
\end{equation*}
\begin{equation}
\gamma_{\textbf{TNB}}=\left(
\begin{array}{c}
1,\cos \left( \int \tau (s)\,ds\right) +a\int \cos \left( \int \tau
(s)\,ds\right) \,ds \\
-\sin \left( \int \tau (s)\,ds\right) ,\cos \left( \int \tau (s)\,ds\right) +
\\
a\int \sin \left( \int \tau (s)\,ds\right) ds+\sin \left( \int \tau
(s)\,ds\right)
\end{array}
\right).
\end{equation}

\subsection{Smarandache curves of Anti-Salkowski curve}

Let $\delta (s)$ be Anti-Salkowski curve in $G_{3}$ with ($\kappa=\kappa (s), \tau
=a=const.$) which can be written as
\begin{equation}
\delta (s)=\left(
\begin{array}{c}
s,\int \left( \int \kappa (s)\cos (as)ds\right) \,ds, \\
\int \left( \int \kappa (s)\sin (as)ds\right) \,ds
\end{array}
\right).
\end{equation}
It gives us the following derivatives
\begin{equation*}
\delta ^{\prime }(s)=\left( 1,\int \kappa (s)\cos (as)ds,\int \kappa (s)\sin
(as)ds\right),
\end{equation*}
\begin{equation*}
\delta ^{\prime \prime }(s)=\left( 0,\kappa (s)\cos (as),\kappa (s)\sin
(as)\right),
\end{equation*}
\begin{equation}
\delta ^{\prime \prime \prime }(s)=\left(
\begin{array}{c}
0,\kappa ^{\prime }(s)\cos (as)-a~\kappa (s)\sin (as), \\
\kappa ^{\prime }(s)\sin (as)+a~\kappa (s)\cos (as)%
\end{array}
\right).
\end{equation}
Further, we obtain the following Frenet vectors $\textbf{T}$, $\textbf{N}$, $\textbf{B}$ in the form
\begin{equation*}
\textbf{T}_{\delta}=\left( 1,\int \kappa (s)\cos (as)ds,\int \kappa (s)\sin (as)ds\right),
\end{equation*}
\begin{equation*}
\textbf{N}_{\delta}=\left( 0,\cos (as),\sin (as)\right),
\end{equation*}
\begin{equation}
\textbf{B}_{\delta}=\left( 0,-\sin (as),\cos (as)\right).
\end{equation}
Thus the above computations of Frenet vectors are give Smarandache curves by
\begin{equation*}
\delta_{\textbf{TN}}=\left(
1,\cos (as)+\int \kappa (s)\cos (as)\,ds,
\int \kappa (s)\sin (as)\,ds+\sin (as)\,
\right)
\end{equation*}
\begin{equation*}
\delta_{\textbf{TB}}=\left(
1,\int \kappa (s)\cos (as)\,\,ds-\sin (as),
\cos (as)+\int \kappa (s)\sin (as)ds
\right)
\end{equation*}
\begin{equation}
\delta_{\textbf{TNB}}=\left(
\begin{array}{c}
1,\cos (as)+\int \kappa (s)\cos (as)\,ds-\sin (as), \\
\cos (as)+\int \kappa (s)\sin (as)ds+\sin (as)%
\end{array}
\right)
\end{equation}

\section{Examples}
\begin{example}
Let $\alpha :I\longrightarrow G_{3}$ be an admissible curve and $%
\kappa \neq 0$ of class $C^{2},\tau \neq 0$ of calss $C^{1}$ its curvature
and torsion, respectively written as
\begin{equation*}
\alpha (s)=\left(
s,\frac{s}{10}\left( -2\cos (2\ln s)+\sin (2\ln s)\right) ,
-\frac{s}{10}\left( \cos (2\ln s)+2\sin (2\ln s)\right)
\right)
\end{equation*}
By differentiation, we get
\begin{equation*}
\alpha ^{\prime }(s)=\left( 1,\cos (\ln s)+\sin (\ln s),-\frac{1}{2}\cos
(2\ln s)\right),
\end{equation*}
\begin{equation*}
\alpha ^{\prime \prime }(s)=\left( 0,\frac{\cos (2\ln s)}{s},\frac{\sin
(2\ln s)}{s}\right),
\end{equation*}
\begin{equation*}
\alpha ^{\prime \prime \prime }(s)=\left(
0,-\frac{\cos (2\ln s)+2\sin (2\ln s)}{s^{2}},
\frac{2\cos (2\ln s)-\sin (2\ln s)}{s^{2}}
\right).
\end{equation*}
Using $(2.5)$ to obtain
\begin{equation*}
\textbf{T}_{\alpha}=\left( 1,\cos (\ln s)\sin (\ln s),-\frac{1}{2}\cos (2\ln s)\right),
\end{equation*}
\begin{equation*}
\textbf{N}_{\alpha}=\left( 0,\cos (2\ln s),\sin (2\ln s)\right),
\end{equation*}
\begin{equation*}
\textbf{B}_{\alpha}=\left( 0,-\sin (2\ln s),\cos (2\ln s)\right).
\end{equation*}
The natural equations of this curve are given by
\begin{equation*}
\kappa_{\alpha}=\frac{1}{s}, \tau_{\alpha}=\frac{2}{s}.
\end{equation*}
Thus, the Smarandache curves of $\alpha$ are respectively, given by
\begin{equation*}
\alpha_{\textbf{TN}}=\left(
1,\cos (2\ln s)+\cos (\ln s)\sin (\ln s),
-\frac{1}{2}\cos (2\ln s)+\sin (2\ln s)
\right),
\end{equation*}
\begin{equation*}
\alpha_{\textbf{TB}}=\left( 1,-\cos (\ln s)\sin (\ln s),\frac{1}{2}\cos (2\ln s)\right),
\end{equation*}
\begin{equation*}
\alpha_{\textbf{TNB}}=\left(
1,\cos (2\ln s)-\cos (\ln s)\sin (\ln s),
\frac{1}{2}\cos (2\ln s)+\sin (2\ln s)
\right).
\end{equation*}
The curve $\alpha$ and their Smarandache curves are shown in Figures 1,2.
\end{example}
\begin{center}
\begin{figure}[!h]
\centering
\includegraphics[scale=0.7]{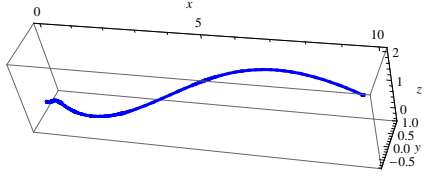}
\label{fig:gh}
\caption{The general helix $\alpha$ in $G_{3}$ with $\kappa =\frac{1}{s}$ and $\tau =\frac{2}{s}$.}
\end{figure}
\end{center}
\begin{center}
\begin{figure}[!h]
\begin{minipage}[t]{0.3\textwidth}
\hspace{-0.1\textwidth}
\centering
\includegraphics[scale=0.7]{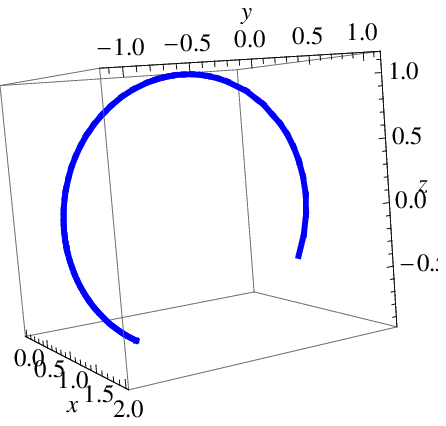}
\label{fig:tngh}
\end{minipage}
\begin{minipage}[t]{0.3\textwidth}
\hspace{-0.1\textwidth}
\centering
\includegraphics[scale=0.7]{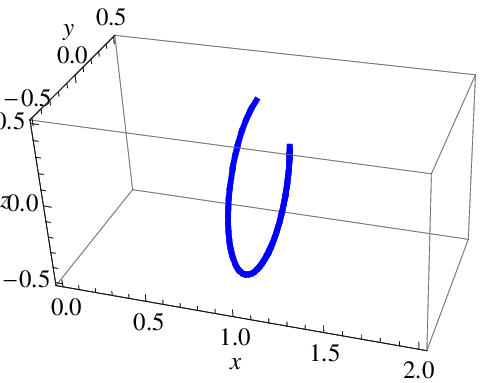}
\label{fig:tbgh}
\end{minipage}
\begin{minipage}[t]{0.3\textwidth}
\hspace{-0.1\textwidth}
\centering
\includegraphics[scale=0.7]{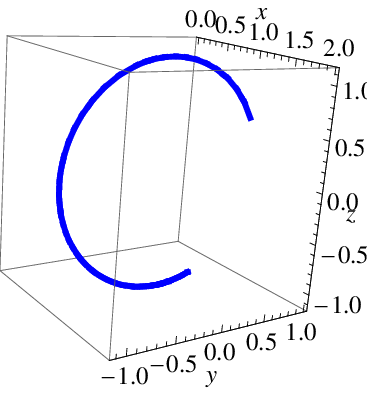}
\label{fig:tnbgh}
\end{minipage}
\caption{From left to right, the $\textbf{TN}$, $\textbf{TB}$ and $\textbf{TNB}$ Smarandache curves of $\alpha$.}
\end{figure}
\end{center}
\newpage
\begin{example}
For an admissible curve $\delta (s)$ in $G_{3}$ parameterized by
\begin{equation*}
\delta (s)=\left(
s,\frac{e^{-s}}{25}(-3\cos (2s)-4\sin (2s)),
\frac{e^{-s}}{25}(4\cos (2s)-3\sin (2s))
\right),
\end{equation*}
we use the derivatives of $\delta ; \delta ^{\prime }, \delta ^{\prime \prime }, \delta ^{\prime \prime \prime }$ to get the associated trihedron of $\delta$ as follows
\begin{equation*}
\textbf{T}_{\delta}=\left\{
1,-\frac{e^{-s}}{5}(\cos (2s)-2\sin (2s)),
-\frac{e^{-s}}{5}(2\cos (2s)+\sin (2s))
\right\},
\end{equation*}
\begin{equation*}
\textbf{N}_{\delta}=\left( 0,\cos (2s),\sin (2s)\right),
\end{equation*}
\begin{equation*}
\textbf{B}_{\delta}=\left( 0,-\sin (2s),\cos (2s)\right).
\end{equation*}
Curvature $\kappa(s)$ and torsion $\tau(s)$ are obtained as follows
\begin{equation*}
\kappa_{\delta}=e^{-s}, \tau_{\delta}=2.
\end{equation*}
According to the above calculations, Smarandache curves of $\delta $ are
\begin{equation*}
\delta_{\textbf{TN}}=\left(
\begin{array}{c}
1,\cos (2s)-\frac{1}{5}e^{-s}(\cos (2s)-2\sin (2s)), \\
\sin (2s)-\frac{1}{5}e^{-s}(2\cos (2s)+\sin (2s))
\end{array}
\right),
\end{equation*}
\begin{equation*}
\delta_{\textbf{TB}}=\left(
\begin{array}{c}
1,-\frac{e^{-s}}{5}\left( \cos (2s)+\left( -2+5e^{s}\right) \sin (2s)\right)
, \\
\cos (2s)-\frac{e^{-s}}{5}(2\cos (2s)+\sin (2s))%
\end{array}
\right),
\end{equation*}
\begin{equation*}
\delta_{\textbf{TNB}}=\left(
\begin{array}{c}
1,\cos (2s)-\frac{e^{-s}}{5}(\cos (2s)-2\sin (2s))-\sin (2s), \\
\cos (2s)+\sin (2s)-\frac{e^{-s}}{5}(2\cos (2s)+\sin (2s))
\end{array}
\right).
\end{equation*}
\end{example}

\begin{center}
\begin{figure}[!h]
\centering
\includegraphics[scale=0.7]{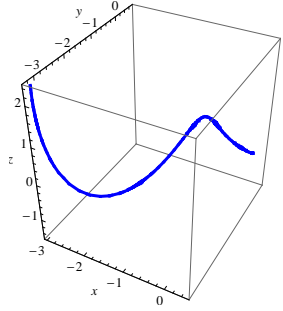}
\label{fig:as}
\caption{The Anti-Salkowski curve $\delta$ in $G_{3}$ with $\kappa_{\delta}=e^{-s}$  and $ \tau_{\delta}=2$.}
\end{figure}
\end{center}
\begin{center}
\begin{figure}[!h]
\begin{minipage}[t]{0.3\textwidth}
\hspace{-0.1\textwidth}
\centering
\includegraphics[scale=0.6]{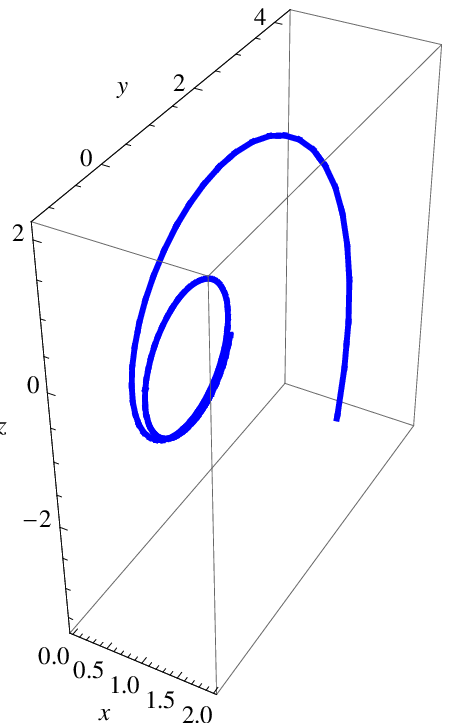}
\label{fig:tnas}
\end{minipage}
\begin{minipage}[t]{0.3\textwidth}
\hspace{-0.1\textwidth}
\centering
\includegraphics[scale=0.6]{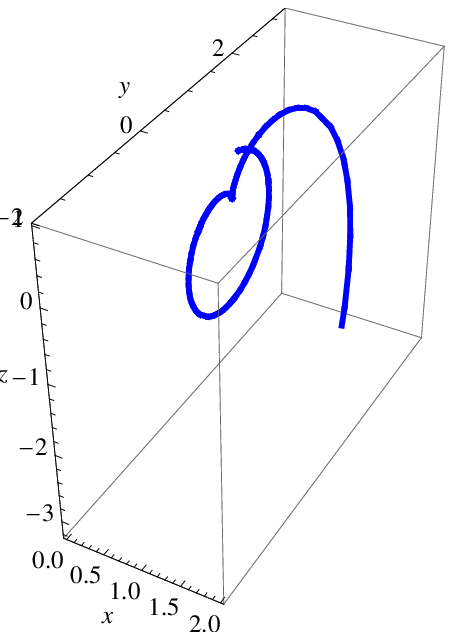}
\label{fig:tbas}
\end{minipage}
\begin{minipage}[t]{0.3\textwidth}
\hspace{-0.1\textwidth}
\centering
\includegraphics[scale=0.6]{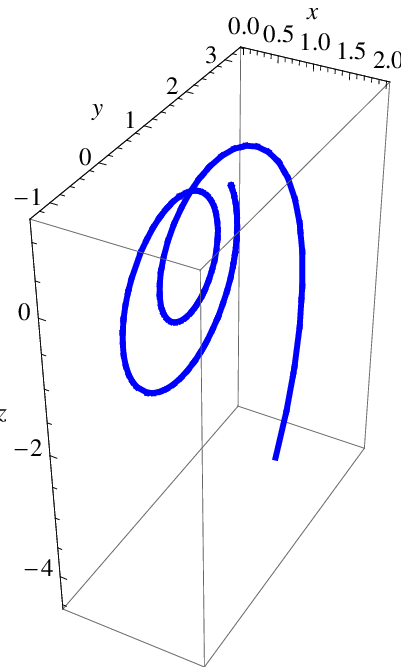}
\label{fig:tnbas}
\end{minipage}
\caption{The $\textbf{TN}$, $\textbf{TB}$ and $\textbf{TNB}$ Smarandache curves of $\delta$.}
\end{figure}
\end{center}

\section{Conclusion}

In the three-dimensional Galilean space, Smarandache curves of an arbitrary curve and some special
curves such as helix, circular helix, Salkowski and Ant-Salkowski curves have been studied. To confirm our main results, two examples (helix and Anti-Salkowski curves) have been given and illustrated.

\end{document}